\documentclass[final]{siamltex}
\usepackage{amsmath,mathrsfs}
\usepackage{amssymb}
\usepackage{graphicx}
\usepackage{hyperref}
\usepackage{pstricks,pst-plot,ifsym}
\usepackage[all]{xy}
\newcommand\subscr[2]{#1_{\textup{#2}}}

\title{Stabilization of systems with one degree of underactuation with energy shaping, a geometric approach\thanks{accepted for publication in SIAM Journal of Control and Optimization}}

\author{Bahman Gharesifard\thanks{Bahman
    Gharesifard is with the Department of
    Mechanical and Aerospace Engineering, University of California San
    Diego, \texttt{bgharesifard@ucsd.edu}}}

\begin{document}

\maketitle

\begin{abstract}
A geometric formulation for stabilization of systems with one degree
of underactuation which fully solves the energy shaping problem for
these system is given. The results show that any linearly
controllable simple mechanical system with one degree of
underactuation is stabilizable by energy shaping, possibly via a
closed-loop metric which is not necessarily positive-definite. An
example of a system with one degree of underactuation is provided
for which the stabilization by energy shaping method is not
achievable using a positive-definite closed-loop metric.
\end{abstract}



\pagestyle{myheadings} \thispagestyle{plain} \markboth{BAHMAN
GHARESIFARD}{STABILIZATION OF SYSTEMS WITH ONE DEGREE OF
UNDERACTUATION}

\section{Introduction}
One of the recent developments for stabilization of simple
mechanical systems is stabilization by energy shaping method. The
central idea concerns the construction of a feedback for which the
closed-loop system inherits the structure of a mechanical system. If
such a feedback exists, the stability of the equilibria can be
guaranteed if the Hessian of the closed-loop potential function is
positive-definite. An important feature of the method, in case it is
applicable, is providing a procedure that allows the design of
nonlinear stabilizing feedbacks.

The first classical appearance of the notion of potential energy
shaping problem is in~\cite{Takegaki_arimoto}. The investigation on
the capabilities of this method continued in the Hamiltonian
framework by looking at the properties of interconnected mechanical
systems. The method is modified into the IDA-PBC method by
introducing the notion of kinetic energy shaping~\cite{Ortega:2002}.
An equivalent version of the IDA-PBC method in Lagrangian
framework, so-called the Controlled Lagrangian method, is initiated
by Bloch, Leonard, Marsden and Chang,~\cite{Bloch1:2000,
Bloch2:2001} and its equivalence to the IDA-PBC has been proved
in~\cite{Chang:2002,Ortega_equ}. In recent work, Chang, Woolsey and
others have realized that the space of possible kinetic energy
feedbacks can be enlarged by considering the addition of appropriate
gyroscopic forcing~\cite{Chang_thesis,Gyr}. It turns out that a
necessary condition for stabilization of simple mechanical control
systems by energy shaping method is linear controllability. For
linear systems, linear controllability is also a sufficient
condition for existence of a stabilizing feedback
\cite{Zenkov_2002_MTNS, Vander_linear}. In
both methods, the question of energy shaping for a mechanical system
reduces to solving a nonlinear system of partial differential
equations. A large number of papers on energy shaping method deals
with finding a parametrization of solutions to this system of
partial differential equations for a particular class of mechanical
systems, for examples
see~\cite{Bloch3,zenkov:2003,Ortega_2001_shaping_revisited,Ortega_Spon_2000}.

A differential geometric approach to the kinetic energy shaping
problem\textemdash the so-called $\lambda$-method\textemdash has
been presented in~\cite{Auckly,Auckly2001,Auckly2002}. A system of
linear partial differential equations is proposed for the kinetic
energy shaping problem in terms of a new variable, $ \lambda
=\subscr{\mathbb{G}}{cl}^{\sharp}\subscr{\mathbb{G}}{ol}^{\flat} $,
where $ \subscr{\mathbb{G}}{ol} $ and $ \subscr{\mathbb{G}}{cl} $
are the open-loop and closed-loop metrics, respectively. The main
idea of the $ \lambda $-method is that it transforms the set of
quasi-linear equations for kinetic energy shaping into a set of
overdetermined linear partial differential
equations~\cite{Auckly2001}. This method is extended to systems with
gyroscopic forces, see~\cite{Chang_thesis}.

Lewis~\cite{Lewis_shaping} has introduced an affine differential
geometric approach to energy shaping. In recent work, sufficient
conditions for the existence of potential energy shaping are derived
assuming that kinetic energy shaping has been
performed~\cite{Lewis_pot}. The results are based on the
integrability theory for linear partial differential equations
developed by Goldschmidt~\cite{Gold_existence} and
Spencer~\cite{Spencer_overestimated}. Gharesifard, et
al.~\cite{Gharesifard:2008} initiated a more systematic geometric
exploration of the kinetic energy shaping partial differential
equations. The authors provide a set of sufficient conditions for
the kinetic energy shaping. Moreover, the kinetic energy shaping
results are coupled with the integrability results of potential
energy shaping~\cite{Lewis_pot} in order to provide a general
approach for the total energy shaping. The technicalities of the
geometric analysis of partial differential equations
in~\cite{Gharesifard:2008} and~\cite{Lewis_pot} might make it hard
to comprehend the importance of these sufficient conditions. The
results of the current paper fully relies on the integrability
analysis of energy shaping partial differential equations and
clarifies the importance of such results.

Numerous systems considered in the literature on energy shaping have
one degree of underactuation. In \cite{Acosta_Ortega_2005_1d} the
authors partially show that, under some conditions, systems with one
degree of underactuation can be stabilized using energy shaping
feedback. The results rely on a transformation of the system of
partial differential equations and do not give any geometric insight
into the energy shaping partial differential equations.
In~\cite{Gharesifard:2008} the authors show that if $
\subscr{\Sigma}{ol} $ is a simple mechanical control system with one
degree of underactuation, for each bundle automorphism which
satisfies the $ \lambda $-equation~\cite{Auckly,Auckly2001}, there
exists a closed-loop metric and a closed-loop potential function
which satisfy the energy shaping system of partial differential
equations. Chang in~\cite{Chang_IFAC_2008} proves this result independently using 
the Cauchy\textendash Kowalevski theorem. Note that this does not guarantee that such solutions are
stabilizing ones. In fact, as we will demonstrate in this paper, the claim of~\cite{Chang_IFAC_2008}  that the existence of such solutions to the kinetic energy shaping guarantees the existence of a positive-definite closed-loop metric which gives rise to an stabilizing energy shaping feedback is not true. However, in the current paper, we show that a linearly controllable simple
mechanical control system with one degree of underactuation can be
stabilized using an energy shaping feedback, possibly via a
closed-loop metric which is not positive-definite.

The paper is organized as follows. In
Section~\ref{section:statement_ene} we recall the affine geometric
formulation of the energy shaping problem. We give a summary of the
integrability results for the partial differential equations in
potential energy shaping in Section~\ref{section:int_result}. We,
directly and without proof, use the results
of~\cite{Gharesifard:2008} and~\cite{Lewis_pot}; thus we do not
review the main integrability theorem of
Goldschmidt~\cite{Gold_existence}. A reader interested in
understanding the formal integrability of partial differential
equations is referred
to~\cite{Gold_existence,Gold_integ,Gold_prog1,Gold_prog2,Spencer_overestimated,Pom2,Guil_alg_inv,Guil_strenberg,seiler_thesis}.
Furthermore, we recall the $ \lambda $-method for kinetic energy
shaping problem. Section~\ref{section:stab_1_du} contains the main
contribution of this paper: we show that all linearly controllable
simple mechanical control systems with one degree of underactuation
are stabilizable using energy shaping method. We fully characterize
a set of solutions to the kinetic energy shaping problem which is
large enough to guarantee the stabilization by energy shaping
method.
\\*
\\*
\textbf{Notation.} The differential geometric notions used in
modeling of simple mechanical systems are assumed here, and the
unfamiliar reader is referred
to~\cite{bullo-lewis,Abraham,foundations_of_mech, Lee,Kobayashi} for
more details. The identity map for a set $ S $ is denoted by $
\mathrm{id}_{S} $ and the image of a map $ f: S \rightarrow W $ by $
\mathsf{Im}(f) $. For a vector space $ \mathsf{V} $ the set of $
(r,s) $-tensors on $ \mathsf{V} $ is denoted by $
\mathrm{T}_s^r(\mathsf{V}) $. By $ \mathsf{S}_k\mathsf{V} $ and $
\Lambda_k\mathsf{V} $ we denote, respectively, the set of symmetric
and skew-symmetric $ (0,k) $-tensors on $ \mathsf{V} $. Let $
\mathsf{A} $ be a $ (0,2) $-tensor on $ \mathsf{V} $. We define the
flat map $ \mathsf{A}^{\flat}: \mathsf{V} \rightarrow \mathsf{V}^* $
by $ \langle\mathsf{A}^{\flat}(u);v\rangle=\mathsf{A}(u,v) $, $
u,v\in \mathsf{V} $. The inverse of the flat map is denoted by $
\mathsf{A}^{\sharp}: \mathsf{V}^* \rightarrow \mathsf{V} $ in case $
\mathsf{A}^{\flat} $ is invertible. We also define a similar
notation for a $ (0,3) $-tensor $ \mathsf{A} $ on $ \mathsf{V} $ by
\[
\langle \mathsf{A}^{\flat}(u),w \rangle =\mathsf{A}(w,u,u), \qquad
u,w \in \mathsf{V}.
\]
For $ \mathsf{S} \subset \mathsf{V} $ and $ \mathsf{W} \subset
\mathsf{V}^* $ we denote
\begin{align*}
\mathrm{ann}(\mathsf{S})&=\{\alpha \in \mathsf{V}^* \ | \ \alpha
(v)=0,\quad
\forall \ v \in \mathsf{S} \},\\
\mathrm{coann}(\mathsf{W})&=\{v \in \mathsf{V} \ | \ \alpha(v)=0,
\quad \forall \ \alpha \in \mathsf{W} \}.
\end{align*}

We denote by $ (\mathsf{E},\pi, \mathsf{Q}) $ a fibered manifold $
\pi : \mathsf{E} \rightarrow \mathsf{Q} $. The \emph{vertical
bundle} of the fibered manifold $ \pi $ is the subbundle of $
\mathsf{T}\pi:\mathsf{T}\mathsf{E}\rightarrow \mathsf{T}\mathsf{Q} $
given by $ \textsf{V}\pi=\ker(\mathsf{T}\pi) $. We denote by $
\mathsf{J}_k\pi $ the \emph{bundle of k-jets} \cite{saunders}. A
local section of $ \pi $ is a pair $ (U, \xi) $, where $ U $ is an
open submanifold of $ \mathsf{Q} $ and $ \xi $ is a map $ \xi: U
\rightarrow \mathsf{E} $ such that $ \pi\circ \xi=\mathrm{id}_{U} $.
If $ (\xi, U) $ is an analytic local section of $ \pi $, we denote
its $ k $-jet by $ j_k \xi $. We denote an element of $ \mathsf{J}_k
\pi $ by $ j_k\xi(x) $, where $ x\in U $. For more information about
geometric properties of jet bundles see~\cite{saunders}. A
\emph{partial differential equation} is a fibered submanifold $
\mathsf{R}_k\subset \mathsf{J}_k\pi $. Goldschmidt theorem
investigates the conditions under which one can construct formal
solutions of a given partial differential equation by
constructing their Taylor series order by
order~\cite{Gold_existence}.

\section{Statement of the energy shaping problem}\label{section:statement_ene} We assume that the reader is familiar with the
affine geometric setup for simple mechanical systems
\cite{bullo-lewis}. A \emph{forced simple mechanical system} is a
quadruple $ \Sigma=(\mathsf{Q}, \mathbb{G}, V,
\subscr{\mathcal{F}}{e}) $, where $ \mathsf{Q} $ is an $ n
$-dimensional manifold called the \emph{configuration manifold}, $
\mathbb{G} $ is a Riemannian metric on $ \mathsf{Q} $, $ V $ is a
function on the configuration manifold called the \emph{potential
function} and $ \subscr{\mathcal{F}}{e}: \mathsf{T}\mathsf{Q}
\rightarrow \mathsf{T}^*\mathsf{Q} $ is a bundle map over $
\mathrm{id}_{\mathsf{Q}} $ called the \emph{external force}. We
denote by $ \nabla^{\mathbb{G}} $ the covariant derivative with
respect to the associated Levi-Civita connection. The governing
equations for a forced simple mechanical system are
\begin{equation*}
\nabla^{\mathbb{G}}_{\gamma'(t)}\gamma'(t)=-\mathbb{G}^{\sharp}\circ
dV(\gamma(t))+\mathbb{G}^{\sharp}\subscr{\mathcal{F}}{e}(\gamma'(t)),
\end{equation*}
where $ \gamma:I\rightarrow \mathsf{Q} $ is an analytic curve on $
\mathsf{Q} $.

Similarly, a \emph{simple mechanical control system} is a quintuple
$ \Sigma=(\mathsf{Q}, \mathbb{G}, V, \subscr{\mathcal{F}}{e},
\mathcal{W}) $ where $ \mathsf{Q} $ is an $ n $-dimensional manifold
called the \emph{configuration manifold}, $ \mathbb{G} $ is a
Riemannian metric on $ \mathsf{Q} $, $ V $ is a function on the
configuration manifold called the \emph{potential function}, $
\subscr{\mathcal{F}}{e}: \mathsf{T}\mathsf{Q} \rightarrow
\mathsf{T}^*\mathsf{Q} $ is a bundle map over $
\mathrm{id}_{\mathsf{Q}} $ called the \emph{external force} and $
\mathcal{W} $ is a subbundle of $ \mathsf{T}^*\mathsf{Q} $ called
the control subbundle~\cite{bullo-lewis}. The governing equations
for a simple mechanical control system are
\begin{equation*}
\nabla^{\mathbb{G}}_{\gamma'(t)}\gamma'(t)=-\mathbb{G}^{\sharp}\circ
dV(\gamma(t))+\mathbb{G}^{\sharp}\subscr{\mathcal{F}}{e}(\gamma'(t))+
\mathbb{G}^{\sharp}u(\gamma'(t)),
\end{equation*}
where $ \gamma:I\rightarrow \mathsf{Q} $ is a curve on $ \mathsf{Q}
$ and $ u: \mathsf{T}\mathsf{Q} \rightarrow \mathcal{W} $ is the
control force.

Given an open-loop simple mechanical control system $
\subscr{\Sigma}{ol}=(\mathsf{Q}, \subscr{\mathbb{G}}{ol},
\subscr{V}{ol}, \subscr{\mathcal{F}}{ol}, \subscr{\mathcal{W}}{ol})
$, we seek a control force such that the closed-loop system is a
forced simple mechanical system $ \subscr{\Sigma}{cl}=(\mathsf{Q},
\subscr{\mathbb{G}}{cl}, \subscr{V}{cl}, \subscr{\mathcal{F}}{cl})
$, possibly with some external force. The reason for seeking this as
the closed-loop system is that the stability analysis of the
equilibria for mechanical systems is well understood \cite[Chapter
6]{bullo-lewis}. In this paper we assume that the open-loop external
force $ \subscr{\mathcal{F}}{ol} $ is zero.

It is well-known that presence of gyroscopic forces enlarges the
space of possible closed-loop metrics
\cite{Chang_thesis,Lewis_shaping} while it does not change the total
energy of the closed-loop system. As we see in next sections,
systems with one degree of underactuation can be stabilized using
energy shaping feedback without gyroscopic forces and thus our
statement of energy shaping problem in this paper does not involve
gyroscopic forces. In following, we present the statement of the
energy shaping problem in the absence of gyroscopic forces.
\begin{definition}[Energy shaping problem in the absence of gyroscopic forces]\label{def:energy_shaping}
Let $ \subscr{\Sigma}{ol}=(\mathsf{Q}, \subscr{\mathbb{G}}{ol},
\subscr{V}{ol}, \subscr{\mathcal{F}}{ol}, \subscr{\mathcal{W}}{ol})
$ be an open-loop simple mechanical control system with $
\subscr{\mathcal{F}}{ol}=0 $. If there exists a bundle map $
\subscr{u}{shp}: \mathsf{T}\mathsf{Q} \rightarrow
\subscr{\mathcal{W}}{ol} $ (called control) with $
\subscr{u}{shp}=-\subscr{u}{kin}-\subscr{u}{pot} $ such that the
closed-loop system is a forced simple mechanical system $
\subscr{\Sigma}{cl}=(\mathsf{Q}, \subscr{\mathbb{G}}{cl},
\subscr{V}{cl}, 0) $ and
\begin{enumerate}
\item \label{eq:kinetic_shap} $ \subscr{\mathbb{G}}{ol}^{\sharp}
\circ
\subscr{u}{kin}(\gamma'(t))=\nabla^{\subscr{\mathbb{G}}{cl}}_{\gamma'(t)}\gamma'(t)
-\nabla^{\subscr{\mathbb{G}}{ol}}_{\gamma'(t)}\gamma'(t), $
\item \label{eq:pot_shap} $
\subscr{u}{pot}(\gamma(t))=\mathbb{G}^{\flat}_{ol}\circ
\mathbb{G}^{\sharp}_{cl}\subscr{dV}{cl}(\gamma(t))-\subscr{dV}{ol}(\gamma(t)),
$
\end{enumerate}
then the control $ \subscr{u}{shp} $ is called an \textit{energy
shaping feedback}.
\end{definition}

Throughout this work, we assume that the equilibrium point $ q_0 \in
\mathsf{Q} $ is a regular point for $ \subscr{\mathcal{W}}{ol} $.
Moreover, we assume that the control codistribution $
\subscr{\mathcal{W}}{ol} $ is \emph{integrable}. This assumption is
common in the literature and many examples fall into this case. The
conditions of Definition~\ref{def:energy_shaping} contain as
unknowns the closed-loop metric $ \subscr{\mathbb{G}}{cl} $. One can
observe that these equations involve the first jet of the unknowns.
One can construct concretely a set of first-order partial
differential equations as necessary and sufficient conditions for
the existence of an energy shaping feedback. Let $ \mathcal{W}
 \subset
\mathsf{T}^*\mathsf{Q} $ be a given subbundle and define the
associated $ \subscr{\mathbb{G}}{ol} $-orthogonal projection map $
P\in \Gamma^{\omega}(\mathsf{T}^*\mathsf{Q} \otimes
\mathsf{T}\mathsf{Q}) $ by
\[
\mathsf{Ker}(P)=\subscr{\mathbb{G}}{ol}^{\sharp}\mathcal{W}.
\]

Note that $ P $ completely prescribes $ \mathcal{W}
 $. We apply $ P
$ to the equation from part~\ref{eq:kinetic_shap} of
Definition~\ref{def:energy_shaping} to arrive at the following
equation:
\[
P(\nabla^{\subscr{\mathbb{G}}{cl}}_{\gamma'(t)}\gamma'(t)
-\nabla^{\subscr{\mathbb{G}}{ol}}_{\gamma'(t)}\gamma'(t))=0.
\]
Assume $ \mathsf{Q} $ is an $ n $-dimensional manifold and $
\mathcal{W}
 $ is an integrable codistribution of dimension $ n-m $.
In adapted local coordinates the kinetic energy shaping partial
differential equation is given by
\begin{equation}\label{eq:pde_kinetic}
P_r^a(\subscr{\mathbb{G}}{cl}^{rl}({\subscr{\mathbb{G}}{cl}}_{,lj,k}+{\subscr{\mathbb{G}}{cl}}_{,lk,j}-{\subscr{\mathbb{G}}{cl}}_{,kj,l})-\subscr{\mathbb{G}}{ol}^{rl}({\subscr{\mathbb{G}}{ol}}_{,lj,k}+{\subscr{\mathbb{G}}{ol}}_{,lk,j}-
{\subscr{\mathbb{G}}{ol}}_{,kj,l}))=0,
\end{equation}
where $ i,j,k,l, r\in\{1,\cdots, n\} $, $ a\in\{1,\cdots, m \} $ and
we denote the first derivative of $ {\subscr{\mathbb{G}}{cl}}_{lj} $
with respect to $ q^k $ by $ {\subscr{\mathbb{G}}{cl}}_{,lj,k} $.
Similarly, let $ \hat{P}: \mathsf{T}^*\mathsf{Q}\rightarrow
\mathsf{T}^*\mathsf{Q}/ \subscr{\mathcal{W}}{ol} $ be the canonical
projection on to the quotient vector bundle. We have
\[
\hat{P}(\subscr{\mathbb{G}}{ol}^{\flat}\circ
\subscr{\mathbb{G}}{cl}^{\sharp}\subscr{dV}{cl}(\gamma(t))-\subscr{dV}{ol}(\gamma(t)))=0.
\]
In local coordinates we have
\begin{equation}\label{eq:pde_pot}
\hat{P}_a^i({\subscr{\mathbb{G}}{ol}}_{,ij}
{\subscr{\mathbb{G}}{cl}}^{jk}{\subscr{V}{cl}}_{,k}-{\subscr{V}{ol}}_{,k})=0,
\end{equation}
where $ i,j,k\in\{1,\cdots, n\} $, $ a\in\{1,\cdots, m \} $ and we
denoted the first derivative of $ {\subscr{V}{cl}} $ with respect to
$ q^k $ by $ {\subscr{V}{cl}}_{,k} $. For more details on the affine
differential geometric setup of energy shaping problem see
\cite{Lewis_shaping}.

\section{Summary of some integrability
results}\label{section:int_result} In~\cite{Lewis_pot} the potential
energy shaping partial differential equation has been shown to be
formally integrable under a surjectivity condition. An important
corollary of this is that the choice of $ \subscr{\mathbb{G}}{cl} $
affects the set of solutions that one might get for potential energy
shaping. A bad choice of $ \mathbb{G}_{cl} $ might make it
impossible to find any potential energy shaping feedback. In a
recent paper \cite{Gharesifard:2008}, the authors show that the
system of partial differential equations for kinetic energy shaping
is formally integrable under a surjectivity condition. Moreover,
they investigate the obstruction for integrability of the total
energy shaping partial differential equations.

Since the integrability conditions of the potential energy shaping
partial differential equations is an integral part of
Theorem~\ref{theorem:1du:main}, we review this result in this
section, without mentioning the proofs. Furthermore, we briefly
recall the $ \lambda $-method for proceeding with the kinetic energy
shaping. We refer an interested reader to~\cite{Gharesifard:2008}
for more details on the integrability character of the partial
differential equations in $ \lambda $-method.

\subsection{Potential energy shaping}
\label{potential_energy_shaping} In this section, we explore aspects
of potential energy shaping. We recall the results for potential
energy shaping after kinetic energy shaping from \cite{Lewis_pot}.
Denote the bundle automorphism $
\subscr{\mathbb{G}}{ol}^{\flat}\circ
\subscr{\mathbb{G}}{cl}^{\sharp} $ by $ \subscr{\Lambda}{cl} $.
Define a codistribution $
\subscr{\mathcal{W}}{cl}=\subscr{\Lambda}{cl}^{-1}(\subscr{\mathcal{W}}{ol})
$ and assume that this codistribution is integrable. Let $
\textsf{PS}\doteq(\mathsf{Q}\times \mathbb{R}, \pi, \mathsf{Q}) $ be
the trivial vector bundle over $ \mathsf{Q} $, so that a section of
$ \pi $ corresponds to a potential function via the formula $ q
\mapsto (q,V(q)) $. We define a $ \mathsf{T}^*\mathsf{Q} $-valued
differential operator $ \mathfrak{D}_{\mathsf{d}}(V)=\mathsf{d}V $
which induces a vector bundle map $ \subscr{\Phi}{pot}:
\mathsf{J}_1\pi \rightarrow \mathsf{T}^*\mathsf{Q} $ such that $
\mathfrak{D}_{\mathsf{d}}(V)(q)=\subscr{\Phi}{pot}(j_1V(q)) $. We
denote by
\[
\pi_{\subscr{\mathcal{W}}{cl}}: \mathsf{T}^*\mathsf{Q}\rightarrow
\mathsf{T}^*\mathsf{Q} / \subscr{\mathcal{W}}{cl}
\]
the canonical projection.
\begin{definition}
Let $ \subscr{\Sigma}{ol}=(\mathsf{Q}, \subscr{\mathbb{G}}{ol},
\subscr{V}{ol}, \subscr{\mathcal{F}}{ol}, \subscr{\mathcal{W}}{ol})
$ be an open-loop simple mechanical control system. The submanifold
$ \subscr{\mathsf{R}}{pot} \subset \mathsf{J}_1\pi $ defined by
\[
\subscr{\mathsf{R}}{pot}=\{p\in \mathsf{J}_1\pi \ | \
\pi_{\subscr{\mathcal{W}}{cl}} \circ
\subscr{\Phi}{pot}(p)=\pi_{\subscr{\mathcal{W}}{cl}}\circ
\Lambda_{cl} ^{-1} \subscr{\mathsf{d}V}{ol}\}
\]
is called the \emph{potential energy shaping submanifold}.
\end{definition}
One can easily observe that the ``equation'' representation of $
\subscr{\mathsf{R}}{pot} $ is given by Equation~(\ref{eq:pde_pot}).

 Let $ \pi_1: \mathsf{J}_1 \pi
\rightarrow \mathsf{Q} $ be the canonical projection. Lewis
\cite{Lewis_pot} gives a set of sufficient conditions under which
the potential shaping problem has a solution. The proof follows from
the integrability theory of partial differential equations; in
particular, the potential energy shaping partial differential
equation has an involutive symbol,
see~\cite{Guil_alg_inv,Gold_integ,Gharesifard:2008} for definition
of involutivity. We recall the definition of $
(\subscr{\mathbb{G}}{ol}$-$\subscr{\mathbb{G}}{cl})$-potential
energy shaping feedback from \cite{Lewis_pot}.
\begin{definition}
A section $ \mathcal{F} $ of $ \mathcal{W} $ is called a \emph{$
(\subscr{\mathbb{G}}{ol}$-$\subscr{\mathbb{G}}{cl})$-potential
energy shaping feedback} if there exists a function $ \subscr{V}{cl}
$ on $ \mathsf{Q} $ such that
\[
\mathcal{F} (q)
=\subscr{\Lambda}{cl}\subscr{\mathsf{d}V}{cl}-\subscr{\mathsf{d}V}{ol},
\qquad q \in \mathsf{Q}.
\]
\end{definition}
The following theorem establishes sufficient conditions for construction of a Taylor series
solution to the potential energy shaping partial differential
equation order-by-order.
\begin{theorem}\label{theorem:Andrew_pot}
Let $ \subscr{\Sigma}{ol}=(\mathsf{Q}, \subscr{\mathbb{G}}{ol},
\subscr{V}{ol}, \subscr{\mathcal{F}}{ol}, \subscr{\mathcal{W}}{ol})
$ be an analytic open-loop simple mechanical control system. Let $
\subscr{\mathbb{G}}{cl} $ be a closed-loop analytic metric. Let $
p_0 \in \subscr{\mathsf{R}}{pot} $ and let $ q_0=\pi_1(p_0) $.
Assume that $ q_0 $ is a regular point for $
\subscr{\mathcal{W}}{ol} $ and that $
\subscr{\mathcal{W}}{cl}=\subscr{\Lambda}{cl}^{-1}\subscr{\mathcal{W}}{ol}
$ is integrable in a neighborhood of $ q_0 $. Then the following
statements are equivalent:
\begin{enumerate}
\item there exists a neighborhood $ U $ of $ q_0 $ and an analytic $
(\subscr{\mathbb{G}}{ol}$-$\subscr{\mathbb{G}}{cl})$ -potential
energy shaping feedback $ \mathcal{F} \in
\Gamma^{\omega}(\mathcal{W}) $ defined on $ U $ which satisfies
\[
\subscr{\Phi}{pot}(p_0)=\subscr{\Lambda}{cl}
\mathsf{d}V(q_0)-\subscr{\mathsf{d}V}{ol}(q_0)+\subscr{\Lambda}{cl}^{-1}
\subscr{\mathsf{d}V}{ol}(q_0),
\]
for a solution $ V $ to $ \subscr{\mathsf{R}}{pot} $;
\item there exists a neighborhood $ U $ of $ q_0 $ such that $ d(\subscr{\Lambda}{cl}^{-1}
\subscr{\mathsf{d}V}{ol})(q)\in\mathsf{I}_2
(\subscr{\mathcal{W}}{cl}|_q) $, where we denote $
\mathsf{I}_2(\subscr{\mathcal{W}}{cl}|_q)=\mathsf{I}(\subscr{\mathcal{W}}{cl}|_q)\cap
\Lambda_2(\mathsf{T}^*_q\mathsf{Q}) $ and the algebraic ideal $
\mathsf{I}(\subscr{\mathcal{W}}{cl}|_q) $ of $
\Lambda(\mathsf{T}^*_q\mathsf{Q}) $ is generated by elements of the
form $ \gamma \wedge\omega $ with $ \gamma \in
\subscr{\mathcal{W}}{cl}|_q $.
\end{enumerate}
\end{theorem}
The theorem gives a set of compatibility conditions for the
existence of a $
(\subscr{\mathbb{G}}{ol}$-$\subscr{\mathbb{G}}{cl})$-potential
energy shaping feedback. Moreover, one can give a full description
of the set of achievable potential energy shaping feedbacks. Let $
\alpha_{cl}=\subscr{\Lambda}{cl}^{-1} \subscr{\mathsf{d}V}{ol} $.
Let us use a coordinate system $ (q^1,\ldots, q^n) $ on $ U $ a
neighborhood of $ q_0 $ such that
\[
\subscr{\mathcal{W}}{cl}|_{q_0} =\mathsf{span}(dq^{m+1}, \cdots,
dq^n).
\]
In these local coordinates we write the one form $ \alpha_{cl} $ as
$ \alpha_{cl}=\alpha_j dq^j $ and compatibility conditions become:
\begin{equation}\label{eq:pot_comp}
\frac{\partial \alpha_j}{\partial q^i}-\frac{\partial
\alpha_i}{\partial q^j}=0 \ \ , \ \ i,j \in \{1, \ldots, m\}.
\end{equation}

\subsection{Kinetic energy shaping (the $ \mathversion{bold}{\lambda}
$-method)} \label{sec:lambda_method}  In following, we recall the
so-called $ \lambda $-method \emph{in the absence of gyroscopic
forces}. The idea is to transform the kinetic energy shaping partial
differential equations to an overdetermined
linear partial differential equation, so-called the $ \lambda
$-equation~\cite{Auckly,Chang_thesis,Gharesifard:2008}.
\begin{theorem}\label{theorem:main_lambda} Let $
\subscr{\Sigma}{ol}=(\mathsf{Q}, \subscr{\mathbb{G}}{ol},
\subscr{V}{ol}, \subscr{\mathcal{F}}{ol}, \subscr{\mathcal{W}}{ol})
$ be an open-loop simple mechanical control system. Let $ P\in
\Gamma^{\omega}(\mathsf{T}^*\mathsf{Q}\otimes \mathsf{T}\mathsf{Q})
$ be the $ \subscr{\mathbb{G}}{ol} $-orthogonal projection as above.
Let $ \subscr{\mathbb{G}}{cl} \in
\Gamma^{\omega}(\mathsf{S}_2^+\mathsf{T}^*\mathsf{Q}) $. If $
\subscr{\mathbb{G}}{ol}^{\flat}=\subscr{\mathbb{G}}{cl}^{\flat}\circ\lambda
$ for $ \lambda \in \Gamma^{\omega}(\mathsf{T}^*\mathsf{Q}\otimes
\mathsf{T}\mathsf{Q}) $, the following two conditions are
equivalent:
\begin{enumerate}
\item $ P(\nabla^{\subscr{\mathbb{G}}{cl}}_{X}X
-\nabla^{\subscr{\mathbb{G}}{ol}}_{X}X)=0 \ \ , \ \ \forall X \in
\Gamma^{\omega}(\mathsf{T}\mathsf{Q}) $;
\item \label{enum:lambda_ab}
\begin{enumerate}
\item $ \nabla^{\subscr{\mathbb{G}}{ol}}_Z(\subscr{\mathbb{G}}{ol}
\lambda)(PX,PY)=0, $ and
\item $ \nabla^{\subscr{\mathbb{G}}{ol}}_{\lambda
PX}\subscr{\mathbb{G}}{cl}(Z,Z)+2\subscr{\mathbb{G}}{cl}(\nabla^{\subscr{\mathbb{G}}{ol}}_Z\lambda
PX,Z)=2\subscr{\mathbb{G}}{ol}(\nabla^{\subscr{\mathbb{G}}{ol}}_{Z}PX,Z)$,
\end{enumerate}
where $ X,Y,Z \in \Gamma^{\omega}(\mathsf{T}\mathsf{Q}) $.
\end{enumerate}
\end{theorem}
For a complete version of the theorem and the proof, in the presence
of gyroscopic forces, see~\cite{Chang_thesis,Gharesifard:2008}. The
set of $ \lambda $-equations have been proved to be formally
integrable under a surjectivity condition~\cite{Gharesifard:2008}.

\subsection{An important corollary for systems with one degree of
underactuaion} For systems with one degree of underactuation the
potential energy shaping partial differential equations is always
formally integrable. The main idea of the proof is that
Equation~(\ref{eq:pot_comp}) vanishes for $ m=1 $, for details of
the proof see \cite{Gharesifard:2008}.

\begin{theorem}\label{theorem:formal_int_1DU}
If $ \subscr{\Sigma}{ol} $ is a simple mechanical control system
with one degree of underactuation, for each bundle automorphism that
satisfies the $ \lambda $-equation, there exists a closed-loop
metric and a closed-loop potential function that satisfy the energy
shaping partial differential equations.
\end{theorem}

In the rest of these paper, we focus on stabilization of the
closed-loop system. Basically, we seek a solution to the energy
shaping partial differential equation for which the Hessian of the
closed-loop potential function can be guaranteed to be
positive-definite.

\section{Stabilization of systems with one degree of
underactuation}\label{section:stab_1_du} In this section, we wish to
determine the stabilizing solutions to the energy shaping partial
differential equations for systems with one degree of
underactuation. Throughout this section, let $ \mathsf{Q} $ be an $
n $-dimensional analytic manifold and $ \subscr{\Sigma}{ol}
=(\mathsf{Q}, \subscr{\mathbb{G}}{}, \subscr{V}{ol},
\subscr{\mathcal{W}}{ol}) $ be an open-loop simple mechanical
control system with one degree of underactuation. We denote the
Hessian of a potential function $ V $ at $ q_0\in \mathsf{Q} $ by $
\mathsf{Hess}(V)(q_0) \in \mathsf{S}_2\mathsf{T}^*_{q_0}\mathsf{Q}
$. In particular, we denote the Hessian of the open-loop potential
function and the closed-loop potential function at the equilibrium
point $ q_0 $ by $ \mathsf{Hess}(\subscr{V}{ol})(q_0) $ and $
\mathsf{Hess}(\subscr{V}{cl})(q_0) $, respectively.

Since the compatibility conditions of
Theorem~\ref{theorem:Andrew_pot} are always satisfied for systems with one
degree of underactuation, Theorem~\ref{theorem:formal_int_1DU}, one
can study the prolongation of the potential energy shaping partial
differential equations instead of the original partial differential
equations. Let $ (q^1,\ldots, q^n) $ be local coordinates in a
neighborhood $ U $ of $ q_0\in \mathsf{Q} $ such that $
\subscr{\mathcal{W}}{ol}=\textrm{span}\{dq^2,\ldots, dq^n\} $ and
let $ P $ be the projection of $ \mathsf{T}^*\mathsf{Q} $ onto $
\mathsf{span}\{ dq^1\} $.

If we prolong the potential energy shaping partial differential
equation and evaluate the result at the origin, noting that $
\mathsf{d}\subscr{V}{cl}(q_0)=0 $, we have
\[
P\big(\subscr{\mathbb{G}}{}^{\flat}(q_0)\subscr{\mathbb{G}}{cl}^{\sharp}(q_0)d^2\subscr{V}{cl}(v)(q_0)-d^2\subscr{V}{ol}(v)(q_0)\big)=0,
\]
where $ v \in \mathsf{T}_{q_0}\mathsf{Q} $, i.e.,
\begin{equation}\label{eq:1du:stab_aff}
\subscr{\mathbb{G}}{cl}^{\sharp}(q_0)\mathsf{Hess}^{\flat}(\subscr{V}{cl})(q_0)
-\subscr{\mathbb{G}}{}^{\sharp}(q_0)\mathsf{Hess}^{\flat}(\subscr{V}{ol})(q_0)
=\subscr{\mathbb{G}}{}^{\sharp}(q_0)(u|_{q_0}),
\end{equation}
where $ u:\mathsf{T}\mathsf{Q}\rightarrow \subscr{\mathcal{W}}{ol}
$. If the system is linearly controllable, then one can design a
control such that $
\subscr{\mathbb{G}}{}^{\sharp}(q_0)\mathsf{Hess}^{\flat}(\subscr{V}{ol})(q_0)
+\subscr{\mathbb{G}}{}^{\sharp}(q_0)(u|_{q_0}) $ is diagonalizable
and positive-definite. It is important to note that this does not
necessarily imply that there exist $ \subscr{\mathbb{G}}{cl} $ and $
\subscr{V}{cl} $ such that $
\subscr{\mathbb{G}}{cl}^{\sharp}(q_0)\mathsf{Hess}^{\flat}(\subscr{V}{cl})(q_0)
$ is positive-definite, since the kinetic energy shaping partial
differential equation puts restrictions on the achievable
closed-loop metrics. However, we will show that, for systems with
one degree of underactuation, the space of solutions of the kinetic
energy shaping partial differential equations is large enough so
that $
\subscr{\mathbb{G}}{cl}^{\sharp}(q_0)\mathsf{Hess}^{\flat}(\subscr{V}{cl})(q_0)
$ can be made positive-definite. We do this in the following steps.
\begin{enumerate}
\item We first identify a simple class of solutions to the $
\lambda $-equation using Proposition~\ref{prop:1DU_lambda_sol}.
\item We show that this class of solutions is large enough to ensure
that Equation~(\ref{eq:1du:stab_aff}) holds with $
\subscr{\mathbb{G}}{cl}^{\sharp}(q_0)\mathsf{Hess}^{\flat}(\subscr{V}{cl})(q_0)
$ diagonalizable and positive-definite.
\end{enumerate}
Let $ U $ be a neighborhood of the equilibrium point $ q_0 \in
\mathsf{Q} $ and let $ (q^1,\ldots, q^n) $ be local coordinates on $
U $. In order to find the class of solutions mentioned in~1, we need
to make some observations about the kinetic energy shaping partial
differential equations for systems with one degree of
underactuation. For these systems, the $ \lambda $-equation in the
adapted local coordinate is given by
\begin{equation}\label{eq:kinetic_1DU}
\frac{\partial}{\partial q^k}
(\mathbb{G}_{1i}\lambda_1^i)-2\mathcal{S}_{k1}^s\mathbb{G}_{si}\lambda_1^i=0,
\end{equation}
where $ \mathcal{S}_{jk}^i $, for $ i,j,k\in\{1,\ldots,n\} $, are
the Levi-Civita connection coefficients associated to $ \mathbb{G} $
and $ i,k,s\in \{1,\ldots, n\} $. Suppose we are seeking solutions
to the $ \lambda $-equation that in local coordinates look like $
\lambda(q)=\lambda_i^jdq^i\otimes \frac{\partial}{\partial q^j} $,
where $ \lambda_i^j\in \mathbb{R} $ and $ q\in U $, i.e., $ \lambda
$ is constant. Then one can write Equation~(\ref{eq:kinetic_1DU}) as
follows:
\begin{equation}\label{eq:kinetic_1DU_spanned}
\left(\frac{\partial \mathbb{G}_{11}}{\partial
q^k}-2\mathcal{S}_{k1}^i\mathbb{G}_{i1}\right)\lambda_1^1+\left(\frac{\partial
\mathbb{G}_{12}}{\partial
q^k}-2\mathcal{S}_{k1}^i\mathbb{G}_{i2}\right)\lambda_1^2+\cdots+\left(\frac{\partial
\mathbb{G}_{1n}}{\partial
q^k}-2\mathcal{S}_{k1}^i\mathbb{G}_{in}\right)\lambda_1^n=0.
\end{equation}
Because $ \mathcal{S} $ is the Levi-Civita connection for $
\mathbb{G} $, the first term vanishes, leaving $ \lambda_1^1 $
arbitrary. One can rewrite Equation~(\ref{eq:kinetic_1DU_spanned})
in the following fashion:
\begin{equation}\label{eq:constant_kientic_cond}
\sum_{i=1}^n\sum_{j=2}^n(\mathcal{S}_{kj}^i\mathbb{G}_{i1}-\mathcal{S}_{k1}^i\mathbb{G}_{ij})\lambda_1^j=0,
\end{equation}
where $ k\in\{1,\ldots,n\} $. Thus, if $ \lambda_2^j=0 $ for $
j\in\{2,\ldots, n\} $, $ \lambda(q) $ is a solution to the $ \lambda
$-equation. Note that we further require that $ \lambda(q)\circ
\mathbb{G}^{\sharp}(q) $ is symmetric. In the following, we describe
the space of such solutions of the $ \lambda $-equation in an
algebraic fashion.

Let $ \mathsf{V} $ be an $ n $-dimensional $ \mathbb{R} $-vector
space and let $ \mathsf{G} \in \mathsf{S}_2\mathsf{V} $ be a
nondegenerate symmetric tensor. Let $ \Phi_{\mathsf{G}}:
\mathsf{V}^*\otimes \mathsf{V} \rightarrow \Lambda_2\mathsf{V} $ be
the map defined by
\[
\Phi_{\mathsf{G}}(\mathsf{A})(v_1,v_2)=\mathsf{A}\circ
\mathsf{G}(v_1,v_2)-\mathsf{A}\circ \mathsf{G}(v_2,v_1),
\]
where $ v_1,v_2\in \mathsf{V} $. The space of all tensors, $
\mathsf{A} \in \mathsf{V}^*\otimes \mathsf{V} $, such that $
\mathsf{A}\circ \mathsf{G} $ is symmetric belongs to the kernel of $
\Phi_{\mathsf{G}} $ and thus is of dimension $ \tfrac{n(n+1)}{2} $,
we denote this subspace by $ \mathsf{S}_{\mathsf{G}}$.\\
Let $ \{e_i\}_{i=1}^n $ be a basis for $ \mathsf{V} $ and let $
\{e^i\}_{i=1}^n $ be its dual. Let $ \mathsf{W}\subset \mathsf{V}^*
$ be the vector subspace generated by $ \{e^2,\ldots, e^n\} $ and
denote its complement by $ \mathsf{E} $. We denote by $
\tilde{\mathsf{S}} $ the space of all $ \mathsf{A}\in
\mathsf{V}^*\otimes \mathsf{V} $ such that, if $ v \in
\mathrm{coann}(\mathsf{W}) $, then $ \mathsf{A}(v) \in
\mathrm{coann}(\mathsf{W}) $, for all $ v \in \mathsf{V} $. A tensor
$ \mathsf{A}\in \tilde{\mathsf{S}} $ can be written as
\[
\mathsf{A}=\mathsf{A}_1^1e^1\otimes
e_1+\sum_{i=2}^n\sum_{j=1}^n\mathsf{A}_i^je^i\otimes e_j,
\]
where $ \mathsf{A}_1^1\in \mathbb{R} $ and $ \mathsf{A}_i^j\in
\mathbb{R} $ for $ i\in \{2,\ldots, n\} $ and $ j \in \{1,\ldots,n\}
$. Thus the dimension of $ \tilde{\mathsf{S}} $ is $ n(n-1)+1 $. If
we denote the restriction of the map $ \Phi_{\mathsf{G}} $ to $
\tilde{\mathsf{S}} $ by $ \Phi_{\mathsf{G}}|_{\tilde{\mathsf{S}}}:
\tilde{\mathsf{S}} \rightarrow \Lambda_2\mathsf{V} $, then $
\ker(\Phi_{\mathsf{G}}|_{\tilde{\mathsf{S}}}) $ is of dimension $
\tfrac{n(n-1)}{2}+1 $. If we additionally require that $
\mathsf{A}\in \ker(\Phi_{\mathsf{G}}|_{\tilde{\mathsf{S}}}) $ be
nondegenerate, we obtain a  $ \tfrac{n(n-1)}{2} $-dimensional
subspace of $ \mathsf{V}^*\otimes \mathsf{V}$.

Let $ \mathsf{Q} $ be an $ n $-dimensional analytic manifold and $
\subscr{\Sigma}{ol} =(\mathsf{Q}, \subscr{\mathbb{G}}{},
\subscr{V}{ol}, \subscr{\mathcal{W}}{ol}) $ be an open-loop simple
mechanical control system with one degree of underactuation. Let $ U
$ be a neighborhood of the equilibrium point $ q_0 \in \mathsf{Q} $
and let $ (q^1,\ldots, q^n) $ be local coordinates on $ U $ such
that $ \subscr{\mathcal{W}}{ol}|_{q}=\mathsf{span}\{dq^2,\ldots,
dq^n\} $, where $ q \in U $. In following, we define a subspace of $
\mathsf{T}_q^*\mathsf{Q}\otimes \mathsf{T}_q\mathsf{Q} $ which is
large enough for stabilization of systems with one degree of
underactuation. Consider the space of solutions to the $ \lambda
$-equation that in local coordinates look like $
\lambda(q)=\lambda_i^jdq^i\otimes \frac{\partial}{\partial q^j} \in
\mathsf{T}^*_{q}\mathsf{Q}\otimes \mathsf{T}_{q}\mathsf{Q} $, where
$ \lambda_i^j\in \mathbb{R} $ and $ q\in U $, and satisfies the
followings
\begin{enumerate}
\item $ \lambda(q)\circ\mathbb{G}^{\sharp}(q) $ is symmetric and nondegenerate;
\item if $ v \in \mathrm{coann}(\mathsf{span}\{ dq^1 \}) $ then $
\lambda(v) \in \mathrm{coann}(\mathsf{span}\{ dq^1 \}) $ for all $ v
\in \mathsf{T}_q\mathsf{Q} $.
\end{enumerate}
We denote this subspace by $ \mathscr{S} $. The following
proposition is a corollary of the algebraic discussion above.
\begin{proposition}\label{prop:1DU_lambda_sol}
$ \mathscr{S} $ is an $ \tfrac{n(n-1)}{2} $-dimensional subspace of
$ \mathsf{T}_q^*\mathsf{Q}\otimes \mathsf{T}_q\mathsf{Q} $.
\end{proposition}

We wish to show that the space of solutions of the $ \lambda
$-equation, described in Proposition~\ref{prop:1DU_lambda_sol}, are
large enough to guarantee that $
\subscr{\mathbb{G}}{cl}^{\sharp}(q_0)\mathsf{Hess}^{\flat}(\subscr{V}{cl})(q_0)
$ can be made diagonalizable and with positive real eigenvalues. If
$ \lambda(q) \in \mathscr{S} $, then
Equation~(\ref{eq:1du:stab_aff}) gives
\begin{align}
&\mathsf{Hess}^{\flat}(\subscr{V}{cl})(q_0)(\frac{\partial}{\partial q^1},\frac{\partial}{\partial q^j})=\frac{1}{\lambda_1^1}\mathsf{Hess}^{\flat}(\subscr{V}{ol})(q_0)(\frac{\partial}{\partial q^1},\frac{\partial}{\partial q^j}),\label{eq:H_1d}\\
&
\subscr{\mathbb{G}}{cl}^{\sharp}(q_0)(dq^1,dq^j)=\lambda_1^1\subscr{\mathbb{G}}{}^{\sharp}(q_0)(dq^1,dq^j)\label{eq:G_1d},
\end{align}
where $ j\in \{1,\ldots,n\} $. As a result, we have the following
proposition.

\begin{proposition}\label{prop:linear_sol_space_aff}
Let $ \mathsf{Q} $ be an $ n $-dimensional analytic manifold and $
\subscr{\Sigma}{ol} =(\mathsf{Q}, \subscr{\mathbb{G}}{},
\subscr{V}{ol}, \subscr{\mathcal{W}}{ol}) $ be an open-loop simple
mechanical control system with one degree of underactuation. Let $ U
$ be a neighborhood of the equilibrium point $ q_0 \in \mathsf{Q} $
and let $ (q^1,\ldots, q^n) $ be local coordinates on $ U $ such
that $ \subscr{\mathcal{W}}{ol}|_{q}=\mathsf{span}\{dq^2,\ldots,
dq^n\} $, where $ q \in U $. Suppose that
\[A=
\subscr{\mathbb{G}}{}^{\sharp}(q_0)\mathsf{Hess}^{\flat}(\subscr{V}{ol})(q_0)
+\subscr{\mathbb{G}}{}^{\sharp}(q_0)(u|_{q_0})
\]
is diagonalizable with real eigenvalues, where $ u|_{q_0}:
\mathsf{T}_{q_0}\mathsf{Q}\rightarrow
\subscr{\mathcal{W}}{ol}|_{q_0} $. Then there exists a closed-loop
metric $ \subscr{\mathbb{G}}{cl} $ and a potential function $
\subscr{V}{cl} $ such that
\begin{enumerate}
\item $
\subscr{\mathbb{G}}{}^{\flat}= \subscr{\mathbb{G}}{cl}^{\flat}\circ
\lambda $, where $ \lambda \in \mathscr{S} $,
\item $
\subscr{\mathbb{G}}{cl}^{\sharp}(q_0)\mathsf{Hess}^{\flat}(\subscr{V}{cl})(q_0)=A
$.
\end{enumerate}
\end{proposition}
\begin{proof}
We only need to show that if~1 holds, then $ \subscr{\mathbb{G}}{cl}
$ and $ \subscr{V}{cl} $ can be selected so that~2 holds. Using
Equations~(\ref{eq:H_1d}) and~(\ref{eq:G_1d}), we can write $
\subscr{\mathbb{G}}{cl}^{\sharp}(q_0) $ in coordinates as
\[
\left(\begin{array}{cc} \lambda_1^1a & \lambda_1^1\mathsf{B} \\
\lambda_1^1\mathsf{B}^{T} & \mathsf{C}
\end{array}\right),
\]
where $ a \in \mathbb{R} $, $ \mathsf{B}\in
\mathrm{L}(\mathbb{R}^{n-1},\mathbb{R}) $, and $ \mathsf{C}\in
\mathsf{S}_2\mathbb{R}^{n-1} $ are such that $
a=\mathbb{G}^{\sharp}(dq^1,dq^1) $ and $
\mathsf{B}(dq^1,dq^j)=\mathbb{G}^{\sharp}(dq^1,dq^j) $ for all $
j\in\{2,\ldots, n\} $. Similarly, $
\mathsf{Hess}^{\flat}(\subscr{V}{cl})(q_0) $ can be written as
\[
\left(\begin{array}{cc} \frac{1}{\lambda_1^1}k &  \frac{1}{\lambda_1^1}\mathcal{B} \\
 \frac{1}{\lambda_1^1} \mathcal{B}^{T} & \mathcal{C}
\end{array}\right),
\]
where $ k\in \mathbb{R} $, $ \mathcal{B}\in
\mathrm{L}(\mathbb{R}^{n-1},\mathbb{R}) $, and $ \mathcal{C}\in
\mathsf{S}_2\mathbb{R}^{n-1} $ are such that $ k =
\mathsf{Hess}^{\flat}(\subscr{V}{ol})(q_0)(\frac{\partial}{\partial
q^1}, \frac{\partial}{\partial q^1}) $ and $ \mathcal{B}(
\frac{\partial}{\partial q^1},\frac{\partial}{\partial q^j})=
\mathsf{Hess}^{\flat}(\subscr{V}{ol})(q_0)(\frac{\partial}{\partial
q^1},\frac{\partial}{\partial q^j}) $ for all $ j\in\{2,\ldots, n\}
$. Thus we have
\[
\subscr{\mathbb{G}}{cl}^{\sharp}(q_0)\mathsf{Hess}^{\flat}(\subscr{V}{cl})(q_0)
=\subscr{\mathbb{G}}{ol}^{\sharp}(q_0)\mathsf{Hess}^{\flat}(\subscr{V}{0l})(q_0)+
\left(\begin{array}{cc} 0 &  0 \\
L_1 &  L_2
\end{array}\right),
\]
where
\begin{enumerate}
\item $
L_1=k\mathsf{B}^T+\frac{1}{\lambda_1^1}\mathsf{C}\mathcal{B}^{T} \in
\mathrm{L}(\mathbb{R},\mathbb{R}^{n-1}) $ and
\item $ L_2=\lambda_1^1\mathsf{B}\mathcal{B}^T+\mathsf{C}\mathcal{C}^T \in
\mathrm{L}(\mathbb{R}^{n-1}\times \mathbb{R}^{n-1}) $
\end{enumerate}
can be set to any value by appropriate choice of $ \mathsf{C} $ and
$ \mathcal{C} $.
\end{proof}

\begin{theorem}\label{theorem:1du:main}
Let $ \subscr{\Sigma}{ol} =(\mathsf{Q}, \subscr{\mathbb{G}}{ol},
\subscr{V}{ol}, \subscr{\mathcal{W}}{ol}) $ be a linearly
controllable open-loop simple mechanical control system with one
degree of underactuation and with $ q_0 \in \mathsf{Q} $ an
equilibrium point. Then the system is stabilizable at $ q_0 $ using
an energy shaping feedback.
\end{theorem}
\begin{proof}
The integrability of the energy shaping partial differential
equations ensures that formal solutions exist. Furthermore,
Theorem~\ref{theorem:formal_int_1DU} implies that the obstructions
of Theorem~\ref{theorem:Andrew_pot} are satisfied for systems with
one degree of underactuation. If the system is linearly
controllable, then one can design a control such that $
\subscr{\mathbb{G}}{}^{\sharp}(q_0)\mathsf{Hess}^{\flat}(\subscr{V}{ol})(q_0)
+\subscr{\mathbb{G}}{}^{\sharp}(q_0)(u|_{q_0}) $ is diagonalizable
and positive-definite. Proposition~\ref{prop:linear_sol_space_aff}
then guarantees that $ \subscr{\mathbb{G}}{cl} $ can be found such
that it satisfies the kinetic energy shaping partial differential
equations, by choosing $ \lambda \in \mathscr{S} $, and taking
\[
\subscr{\mathbb{G}}{cl}^{\sharp}(q_0)\mathsf{Hess}^{\flat}(\subscr{V}{cl})(q_0)=
\subscr{\mathbb{G}}{}^{\sharp}(q_0)\mathsf{Hess}^{\flat}(\subscr{V}{ol})(q_0)
+\subscr{\mathbb{G}}{}^{\sharp}(q_0)(u|_{q_0})
\]
to be diagonalizable with positive real eigenvalues.
\end{proof}

Note that this proof does not require that the closed-loop metric be
positive-definite and in fact, there are cases for which energy
shaping is not possible with positive-definite closed-loop metrics;
an example of this is presented in Example~\ref{example:1d_1}. The
following proposition clarifies when it is necessary to perform
kinetic energy shaping for systems with one degree of
underactuation.
\begin{proposition}\label{prop:1DU_system_no_kes}
Let $ \mathsf{Q} $ be an $ n $-dimensional manifold and let $
\subscr{\Sigma}{ol}=(\mathsf{Q}, \subscr{\mathbb{G}}{},
\subscr{V}{ol}, \subscr{\mathcal{W}}{ol}) $ be a linearly
controllable simple mechanical system. Let $ U $ be a neighborhood
 of $ q_0\in \mathsf{Q} $ such that $
\subscr{\mathcal{W}}{ol}=\mathsf{span}\{dq^2,\ldots, dq^n\} $. If $
\mathsf{Hess}(\subscr{V}{ol})(\frac{\partial}{\partial
q^1},\frac{\partial}{\partial q^1})>0 $, the system can be
stabilized around its equilibrium point $ q_0 $ without kinetic
energy shaping.
\end{proposition}
\begin{proof}
We shall show that $ \subscr{\Sigma}{ol} $ is stabilizable using an
energy shaping feedback with $
\subscr{\mathbb{G}}{cl}=\subscr{\mathbb{G}}{} $.
Equation~(\ref{eq:1du:stab_aff}) then reads
\[
\mathsf{Hess}^{\flat}(\subscr{V}{cl})(q_0)=
\mathsf{Hess}^{\flat}(\subscr{V}{ol})(q_0)+u|_{q_0},
\]
where $ u $ is a feedback. Note that since $
\mathsf{Hess}(\subscr{V}{cl}) $ is symmetric, it is
positive-definite if and only if all of its principal minors are
positive. The first principal minor of $
\mathsf{Hess}^{\flat}(\subscr{V}{cl}) $ is positive. Then, by linear
controllability, one can choose the controls so that the system is
stabilizable at the equilibrium point $ q_0 $, similar to
Proposition~\ref{prop:linear_sol_space_aff}.
\end{proof}

Next, we present an example of energy shaping for simple mechanical
systems with one degree of underactuation for which the energy
shaping is possible \emph{only} via a closed-loop metric that is not
positive-definite.

\section{Example}\label{example:1d_1} Consider the stabilization
problem for a simple mechanical control system $
\Sigma=(\mathbb{R}^2,\mathbb{G}, \subscr{V}{ol}, 0,
\subscr{\mathsf{W}}{ol}) $ at the origin $ q_0=\mathbf{0}\in
\mathbb{R}^2 $, where
\begin{enumerate}
\item $ \subscr{\mathbb{G}}{}=((q^2)^2+1)dq^1\otimes
dq^2+((q^1)^2+1)dq^2\otimes dq^2 $,
\item $ \subscr{V}{ol}=-(q^1)^2+2q^1q^2+(q^2)^2 $, and
\item $ \subscr{\mathcal{W}}{ol}=\mathsf{span}\{ dq^2\} $.
\end{enumerate}
This system is linearly controllable at the origin. We show that,
for any solution of the $ \lambda $-equation, the constant term in
the Taylor expansion of $ \lambda_1^2 $ is always zero. In order to
show this, we need to modify
Equation~(\ref{eq:constant_kientic_cond}) by adding an extra term,
since $ \lambda $, in a neighborhood of $ q_0 $, is not necessarily
chosen from $ \mathscr{S} $. We have
\begin{equation*}
\sum_{i=1}^n(\mathbb{G}_{1i}\frac{\partial \lambda_1^i}{\partial
q^k}+\sum_{j=2}^n(\mathcal{S}_{kj}^i\mathbb{G}_{i1}-\mathcal{S}_{k1}^i\mathbb{G}_{ij})\lambda_1^j)=0,
\end{equation*}
for all $ k\in\{1,\ldots,n\} $. For this example, by substituting
the nonzero Christoffel symbols, we have
\begin{align}
&((q^2)^2+1)\frac{\partial \lambda_1^1}{\partial
q^1}+2q^2\lambda_1^2=0,\label{eq:ex:modified_1d_1}\\
&((q^2)^2+1)\frac{\partial \lambda_1^1}{\partial
q^2}-2q^1\lambda_1^2=0 \label{eq:ex:modified_2d_1}.
\end{align}
It is clear that $ \lambda_1^1(q_0) $ can be chosen arbitrarily.
Consider formal expressions for $ \lambda_1^2 $ and $ \lambda_1^1 $:
\begin{align*}
&\lambda_1^1=C_{00}+C_{10}q^1+C_{01}q^2+C_{20}(q^1)^2+C_{02}(q^2)^2+C_{11}q^1q^2+\cdots,\\
&\lambda_1^2=D_{00}+D_{10}q^1+D_{01}q^2+D_{20}(q^1)^2+D_{02}(q^2)^2+D_{11}q^1q^2+\cdots,
\end{align*}
where $ C_{ij}, D_{ij} \in \mathbb{R} $ for $ i,j \in
\mathbb{Z}_{\geq0} $. If $ \lambda_1^1 $ and $ \lambda_1^2 $ satisfy
Equations~(\ref{eq:ex:modified_1d_1})
and~(\ref{eq:ex:modified_2d_1}), then $ C_{11}=D_{00}=0 $, i.e., $
\lambda_1^2(q_0)=0 $. Thus the closed-loop metric at the origin has
the form $ \subscr{\mathbb{G}}{cl}(q_0)=\frac{1}{a}dq^1\otimes
dq^1+\frac{1}{c}dq^2\otimes dq^2 $, where $ a,c \in
\mathbb{R}\backslash \{0\} $ and $ \lambda_1^1(q_0)=a $.
Equation~(\ref{eq:H_1d}) implies that
\[
\mathsf{Hess}^{\flat}(\subscr{V}{cl})(q_0)=\left(\begin{array}{cc}\tfrac{-2}{a} & \tfrac{2}{a} \\
\tfrac{2}{a} & k\end{array}\right),
\]
where $ k\in \mathbb{R} $. Thus
\[
\subscr{\mathbb{G}}{cl}^{\sharp}(q_0)\mathsf{Hess}^{\flat}(\subscr{V}{cl})(q_0)=\left(\begin{array}{cc}-2 & 2 \\
\tfrac{2c}{a} & ck \end{array}\right).
\]
It is easy to see that one has to choose $ \tfrac{2c}{a}<0 $ and $
ck>2 $ in order to make $
\subscr{\mathbb{G}}{cl}^{\sharp}(q_0)\mathsf{Hess}^{\flat}(\subscr{V}{cl})(q_0)
$ positive-definite, i.e., none of the achievable closed-loop
metrics is positive-definite. However, one can choose $ a,c,k\in
\mathbb{R} $ so that $
\subscr{\mathbb{G}}{cl}^{\sharp}(q_0)\mathsf{Hess}^{\flat}(\subscr{V}{cl})(q_0)
$ is positive-definite, for example $ a=-\tfrac{191}{100} $, $
c=\tfrac{43}{10} $, and $ k=1 $.

\textbf{Remark.} If we take the open-loop metric given by
\[
\subscr{\mathbb{G}}{}=((q^2)^2+1)dq^1\otimes
dq^2+((q^1)^2+1)dq^2\otimes dq^2+2q^1q^2(dq^1\otimes
dq^2+dq^2\otimes dq^1),
\]
then $ \lambda_1^2(q_0) $ need not be zero and the system can be
shown to be stabilizable by the energy shaping method with a
positive-definite closed-loop metric. This reveals that a slight
change in the structure of the open-loop Levi-Civita connection has
a huge impact on the achievable closed-loop metrics.

\section{Summary}
In this paper, we fully solved the problem of stabilization of
systems with one degree of underactuation. The result completely
relies on the integrability analysis of partial differential
equations involved in energy shaping. We illustrated that all
linearly controllable simple mechanical control systems with one
degree of underactuation can be stabilized using an energy shaping
feedback, with closed-loop metrics which are not necessarily
positive-definite. We also characterized the simple mechanical
systems for which the energy shaping is achievable without kinetic
energy shaping. Finally, we gave an example of a simple mechanical
control system with one degree of underactuation for which there
exists no solution to the energy shaping problem with
positive-definite closed-loop metric. The results give some useful
insight about the structure of kinetic energy shaping.

\section*{Acknowledgments}
The author thanks Drs.\ Andrew Lewis and Abdol-Reza Mansouri of
Department of Mathematics and Statistic of Queen's University for
great suggestions and discussions about the results of this paper.
In particular, the author thanks Dr.\ Andrew Lewis for improving the
proof of Proposition~\ref{prop:linear_sol_space_aff}. The author thanks the reviewers for their useful comments.


\clearpage
\begin{thebibliography}{10}

\bibitem{foundations_of_mech}
{\sc R.~Abraham and J.~E. Marsden}, {\em Foundations of Mechanics},
  Addison\textendash Wesley, 2~ed., 1978.

\bibitem{Abraham}
{\sc R.~Abraham, J.~E. Marsden, and T.~S. Ratiu}, {\em Manifolds, Tensor
  Analysis, and Applications}, no.~75 in Applied Mathematical Sciences,
  Springer\textendash Verlag, 2~ed., 1988.

\bibitem{Acosta_Ortega_2005_1d}
{\sc J.~A. Acosta, R.~Ortega, A.~Astolfi, and A.~D. Mahindrakar}, {\em
  Interconnection and damping assignment passivity-based control of mechanical
  systems with underactation degree one}, IEEE, Transactions on Automatic
  Control, 50 (2005), pp.~1936--1955.

\bibitem{Auckly2001}
{\sc D.~R. Auckly and L.~V. Kapitanski}, {\em Mathematical problems in the
  control of underactuated systems}, CRM Proc. Lecture Notes 27, AMS,
  Providence, RI,  (2001), pp.~29--40.

\bibitem{Auckly2002}
\leavevmode\vrule height 2pt depth -1.6pt width 23pt, {\em On the
  $\lambda$-equations for matching control laws}, SIAM Journal of Control and
  Optimization, 41 (2002), pp.~1372--1388.

\bibitem{Auckly}
{\sc D.~R. Auckly, L.~V. Kapitanski, and W.~White}, {\em Control of nonlinear
  underactuated systems}, Communications on Pure and Applied Mathematics, 53
  (2000), pp.~354--369.

\bibitem{Ortega_equ}
{\sc G.~Blankenstein, R.~Ortega, and A.~J. van~der Schaft}, {\em {T}he matching
  conditions of controlled {L}agrangians and {IDA}-passivity based control},
  International Journal of Control, 75 (2002), pp.~645--665.

\bibitem{Bloch3}
{\sc A.~Bloch, A.~E. Leonard, and J.~E Marsden}, {\em Controlled {L}agrangians
  and the stabilization of euler\textendash poincare mechanical systems},
  International Journal of Robust and nonlinear control, 11 (2001),
  pp.~191--214.

\bibitem{Bloch2:2001}
{\sc A.~M. Bloch, D.~E. Chang, N.~E. Leonard, and J.~E. Marsden}, {\em
  Controlled {L}agrangians and the stabilization of mechanical systems. {II}.
  {P}otential shaping}, IEEE Transactions on Automatic Control, 46 (2001),
  pp.~1556--1571.

\bibitem{Bloch1:2000}
{\sc A.~M. Bloch, N.~E. Leonard, and J.~E. Marsden}, {\em Controlled
  {L}agrangians and the stabilization of mechanical systems. {I}. {T}he first
  matching theorem}, IEEE Transactions on Automatic Control, 45 (2000),
  pp.~2253--2270.

\bibitem{bullo-lewis}
{\sc F.~Bullo and A.~D. Lewis}, {\em Geometric Control of Mechanical Systems:
  {M}odeling, Analysis, and Design for Simple Mechanical Control Systems},
  no.~49 in Texts in Applied Mathematics, Springer\textendash Verlag, 2004.

\bibitem{Chang_thesis}
{\sc D.~E. Chang}, {\em Controlled {L}agrangian and {H}amiltonian systems}, PhD
  thesis, California Institute of Technology, 2002.

\bibitem{Chang_IFAC_2008}
\leavevmode\vrule height 2pt depth -1.6pt width 23pt, {\em Some results on
  stabilizability of controlled lagrangian systems by energy shaping}, Proc.
  IFAC World Congress, Seoul Korea,  (2008).

\bibitem{Chang:2002}
{\sc D.~E. Chang, A.~M. Bloch, N.~E. Leonard, J.~E. Marsden, and C.~A.
  Woolsey}, {\em {T}he equivalence of controlled {L}agrangians and controlled
  {H}amiltonian systems}, ESAIM: Control, Optimization and Calculus of
  Variations, 8 (2002), pp.~393--422.

\bibitem{Gharesifard:2008}
{\sc B.~Gharesifard, A.~D. Lewis, and A.~R. Mansouri}, {\em A geometric
  framework for stabilization by energy shaping: Sufficient conditions for
  existence of solutions}, Special Issue Dedicated to the 70th Birthday of
  Roger W. Brockett, Communications in Systems and Information, 8 (2008),
  pp.~353--398.

\bibitem{Gold_existence}
{\sc H.~L. Goldschmidt}, {\em Existence theorems for analytic linear partial
  differential equations}, Annals of Mathematics. Second Series, 86 (1967),
  pp.~246--270.

\bibitem{Gold_integ}
\leavevmode\vrule height 2pt depth -1.6pt width 23pt, {\em Integrability
  criteria for systems of nonlinear partial differential equations}, Journal of
  Differential Geometry, 1 (1967), pp.~269--307.

\bibitem{Gold_prog1}
\leavevmode\vrule height 2pt depth -1.6pt width 23pt, {\em Prolongation of
  linear partial differential equations. {I}, {A} conjecture of \'{E}lie
  {C}artan}, Annales Scientifiques de l'\'Ecole Normale Sup\'erieure.
  Quatri\`eme S\'erie, 1 (1968), pp.~417--444.

\bibitem{Gold_prog2}
\leavevmode\vrule height 2pt depth -1.6pt width 23pt, {\em Prolongation of
  linear partial differential equations. {II}, {I}nhomogeneous equations},
  Annales Scientifiques de l'\'Ecole Normale Sup\'erieure. Quatri\`eme S\'erie,
  1 (1968), pp.~617--625.

\bibitem{Guil_alg_inv}
{\sc V.~Guillemin and M.~Kuranishi}, {\em Some algebraic results concerning
  involutive subspaces}, American Journal of Mathematics, 90 (1968),
  pp.~1307--1320.

\bibitem{Guil_strenberg}
{\sc W.~Guillemin and S.~Sternberg}, {\em An algebraic model of transitive
  differential geometry}, American Mathematical Society. Bulletin. New Series,
  70 (1964), pp.~16--47.

\bibitem{Kobayashi}
{\sc S.~Kobayashi and K.~Nomizu}, {\em Foundations of Differential Geometry},
  vol.~1, Wiley-Interscience, 1981.

\bibitem{Lee}
{\sc John~M. Lee}, {\em Introduction to Smooth Manifolds}, Springer\textendash
  Verlag, 2002.

\bibitem{Lewis_shaping}
{\sc A.~D. Lewis}, {\em Notes on energy shaping}, 43rd IEEE Conference on
  Decision and Control,  (2004), pp.~4818--4823.

\bibitem{Lewis_pot}
\leavevmode\vrule height 2pt depth -1.6pt width 23pt, {\em Potential energy
  shaping after kinetic energy shaping}, 45th IEEE Conference on Decision and
  Control,  (2006), pp.~3339--3344.

\bibitem{Ortega_Spon_2000}
{\sc R.~Ortega and M.~W. Spong}, {\em Stabilization of underactuated mechanical
  systems via interconnection and damping assignment}, IFAC Workshop on
  Lagrangian and Hamiltonian Methods for Nonlinear Control, Princeton, NJ,
  March 16-18,  (2000).

\bibitem{Ortega:2002}
{\sc R.~Ortega, M.~W. Spong, F.~G´omez-Estern, and G.~Blankenstein}, {\em
  Stabilization of a class of underactuated mechanical systems via
  interconnection and damping assignment}, IEEE Transactions on Automatic
  Control, 47 (2002), pp.~1218--1233.

\bibitem{Ortega_2001_shaping_revisited}
{\sc R.~Ortega, A.~J van~der Schaft, I.~Mareels, and B.~Maschke}, {\em Energy
  shaping control revisited}, Springer\textendash Verlag, Berlin, 2001.

\bibitem{Pom2}
{\sc J.~F. Pommaret}, {\em Partial Differential Equations and Group Theory: New
  Perspectives for Applications}, Mathematics and its Applications, Kluwer
  Academic Publishers, 1994.

\bibitem{Vander_linear}
{\sc Stephen Prajna, Arjan van~der Schaft, and Gjerrit Meinsma}, {\em An {LMI}
  approach to stabilization of linear port-controlled {H}amiltonian systems},
  Systems and Control Letters, 45 (2002), pp.~371--385.

\bibitem{saunders}
{\sc D.~J. Saunders}, {\em {T}he Geometry of Jet Bundles}, Cambridge University
  Press, 1989.

\bibitem{seiler_thesis}
{\sc W.~M. Seiler}, {\em Analysis and application of the formal theory of
  partial differential equations}, PhD thesis, Lancaster University, 1994.

\bibitem{Spencer_overestimated}
{\sc D.~C. Spencer}, {\em Overdetermined systems of linear partial differential
  equations}, American Mathematical Society. Bulletin. New Series, 75 (1967),
  pp.~159--193.

\bibitem{Takegaki_arimoto}
{\sc M.~Takegaki and S.~Arimoto}, {\em A new feedback method for dynamic
  control of manipulators}, American Society of Mechanical Engineers.
  Transactions of the ASME. Series G. Journal of Dynamical Systems and
  Measurement Control, 102 (1981), pp.~119--125.

\bibitem{Gyr}
{\sc C.~Woolsey, C.~K. Reddy, A.~Bloch, D.~E. Chang, N.~E. Leonard, and J.~E.
  Marsden}, {\em Controlled {L}agrangian systems with gyroscopic forcing and
  dissipation}, European Journal of Control, 10 (2004).

\bibitem{Zenkov_2002_MTNS}
{\sc D.~V. Zenkov}, {\em Matching and stabilization of linear mechanical
  systems}, Mathematical Theory of Networks and Systems,  (2002).

\bibitem{zenkov:2003}
{\sc D.~V. Zenkov, A.~M. Bloch, and J.~E. Marsden}, {\em {C}ontrolled
  {L}agrangian methods and tracking of accelerated motions}, Proceedings of the
  IEEE International Conference on Decision and Control,  (2003), pp.~533--538.

\end{thebibliography}
\end{document}